\documentclass[11pt,reqno]{amsart}
\usepackage{amssymb,amscd,amsbsy}
\usepackage{amssymb,amscd,amsbsy,mathrsfs}
\renewcommand{\a}{\alpha}
\renewcommand{\b}{\beta}

\newcommand{\e}{\varepsilon}

\renewcommand{\l}{\lambda}

\newcommand{\s}{\sigma}

\newcommand{\f}{\varphi}

\newcommand{\D}{\Delta}
\renewcommand{\L}{\Lambda}

\newcommand{\h}{{\mathscr H}}

\newcommand{\K}{{\mathscr K}}

\newcommand{\X}{{\mathscr X}}

\newcommand{\W}{{\mathscr W}}

\newcommand{\C}{{\Bbb C}}

\newcommand{\Z}{{\Bbb Z}}

\newcommand{\bs}{\boldsymbol}

\newcommand{\bS}{{\boldsymbol S}}

\newcommand{\rf}[1]{(\ref{#1})}

\newcommand{\df}{\stackrel{\mathrm{def}}{=}}
\newcommand{\dist}{\operatorname{dist}}

\newcommand{\spn}{\operatorname{span}}

\newcommand{\rank}{\operatorname{rank}}
\newcommand{\const}{\operatorname{const}}

\newcommand{\eeq}{\end{equation}}
\newcommand{\beq}{\begin{equation}}
\newcommand{\bay}{\begin{eqnarray}}
\newcommand{\ba}{\begin{align*}}
\newcommand{\ea}{\end{align*}}
\newcommand{\ey}{\end{eqnarray}}
\newcommand{\bey}{\begin{eqnarray*}}
\newcommand{\eey}{\end{eqnarray*}}

\newcommand{\be}{\infty}

\newcommand{\bl}{\blacksquare}

\newcommand{\Pf}{{\bf Proof. }}

\newcommand{\ov}{\overline}

\newtheorem{thm}{\hspace{\parindent}Theorem}[section]

\newtheorem{cor}[thm]{\hspace{\parindent}Corollary}
\newtheorem{lem}[thm]{\hspace{\parindent}Lemma}

\theoremstyle{remark}

\newtheorem*{rem*}{Remark}

\newcommand{\ps}{p_\sharp}
\newcommand{\ee}{\boldsymbol{e}}

\newcommand\fM{\frak M}

\newcommand\mX{\mathcal{X}}
\newcommand\mY{\mathcal{Y}}

\newcommand\tp{\ell^\be\otimes_p\ell^\be}

\newcommand{\Bbbone}{{\rm{1\mathchoice{\kern-0.25em}{\kern-0.25em}{\kern-0.2em}{\kern-0.2em}I}}}

\begin{document}

\newcommand{\vse}{\vspace{.2in}}
\numberwithin{equation}{section}

\title{Schur multipliers of Schatten--von Neumann \\
classes $\boldsymbol{S_p}$}
\author{A.B. Aleksandrov and  V.V. Peller}
\thanks{The research of the first author is partially supported by RFBR grant 17-01-00607.
The publication was prepared with the support of the
RUDN University Program 5-100}
\thanks{Corresponding author: V.V. Peller; email: peller@math.msu.edu}

\begin{abstract}
We study in this paper properties of Schur multipliers of Schatten von Neumann classes 
$\bS_p$. We prove that for $p\le1$, Schur multipliers of $\bS_p$ are necessarily completely bounded. 

We also introduce for $p\le1$ a  scale $\W_p$ of tensor products of $\ell^\be$ and prove that matrices in $\W_p$ are Schur multipliers of $\bS_p$.  We compare this sufficient condition with the sufficient condition of membership in the $p$-tensor product of $\ell^\be$ spaces. 
\end{abstract}

\maketitle

\tableofcontents

\setcounter{section}{0}
\section{\bf Introduction}
\setcounter{equation}{0}
\label{In}

\

In this paper we are going to study matrix Schur multipliers of Schatten--von Neumann classes 
$\bS_p$ and, more generally, Schur multipliers with respect to spectral measures.

Recall that for matrices $A=\{a_{jk}\}_{j,k\ge0}$ and $B=\{b_{jk}\}_{j,k\ge0}$, their {\it Schur--Hadamard product} $A\star B$ is defined by
$$
A\star B\df\{a_{jk}b_{jk}\}_{j,k\ge0}.
$$ 
We are going to use the same notation $A\star B$ in the case when 
$A=\{a_{jk}\}_{j,k\ge0}$ is a matrix with complex entries and $B=\{B_{jk}\}_{j,k\ge0}$
is a matrix with operator entries, i.e., the entries $B_{jk}$ are bounded linear operators on a Hilbert space. In this case
$$
A\star B\df\{a_{jk}B_{jk}\}_{j,k\ge0}
$$ 
is an operator matrix.

\medskip

{\bf Definition.}
For $p\in(0,\be]$, we denote by $\fM_p$ the space of {\it matrix Schur multipliers} of $\bS_p$, i.e., a matrix $A=\{a_{jk}\}_{j,k\ge0}$, by definition, belongs to $\fM_p$
if
$$
A\star B\in\bS_p \quad\mbox{for every scalar matrix}\quad B\quad\mbox{in}\quad\bS_p.
$$
As usual, we say that a matrix belongs to $\bS_p$ if it induces a linear operator on the sequence space $\ell^2$ of class $\bS_p$. Note that for convenience, we mean by 
$\bS_\be$  the class of bounded matrices. 
The norm ($p$-norm for $p<1$) $\|A\|_{\fM_p}$ of a matrix Schur multiplier $A$ is defined by
\bay
\label{fMp}
\|A\|_{\fM_p}\df\sup\big\{\|A\star B\|_{\bS_p}:~\|B\|_{\bS_p}\le1\big\}.
\ey

\medskip

Recall that for $p<1$, a functional $\|\cdot\|$ on a vector space $X$ is called a {\it $p$-norm} if it has the following properties:
$$
\begin{array}{l}
\|\a x\|=|\a|\cdot\|x\|,~\quad x\in X,~\a\in\C;\\[.2cm]
\|x\|=0\quad\mbox{only if }\quad x=0;\\[.2cm]
\|x+y\|^p\le\|x\|^p+\|y\|^p,\quad x,~y\in X.
\end{array}
$$
The $p$-norm induces the following metric on $X$:
$$
\dist(x,y)=\|x-y\|^p.
$$
$X$ is called a {\it $p$-Banach space} if it is a complete space with respect to this metric.

The fact that the right-hand side of \rf{fMp} is finite follows easily from the closed graph theorem (it holds not only for Banach spaces but also for complete metric vector spaces, see \cite{Ba}, Ch. 1, \S\,3).

It can easily be verified that for $p<1$, the space $\fM_p$ is a $p$-Banach space with respect to the $p$-norm $\|\cdot\|_{\fM_p}$.

It is well known (see, e.g., \cite{Be}) that $\fM_p=\fM_{p'}$ for $p\ge1$, where 
$p'\df p/(p-1)$.

\medskip

{\bf Definition.} Let $0<p\le\be$. A numerical matrix $A$ is called a {\it completely bounded Schur multiplier} of $\bS_p$ if
$$
A\star B\in\bS_p \quad\mbox{for every operator matrix}\quad B\quad\mbox{in}\quad\bS_p.
$$

\medskip

Obviously, a completely bounded Schur multiplier of $\bS_p$ is necessarily a Schur multiplier of $\bS_p$. The converse trivially holds for $p=2$. It also holds for $p=1$ and $p=\be$, see \cite{Be}.

The famous problem of whether a Schur multiplier of $\bS_p$ must be completely bounded for each $p$  was posed in \cite{Pis}.

One of the main results of this paper is an affirmative answer to this question in the case $p<1$. This will be established in \S\:\ref{vpoogr}. 
The question for $p\in(1,2)\cup(2,\be)$ remains open.

To prove the main result of \S\:\ref{vpoogr}, we obtain in \S\:\ref{dois} a general result 
on double operator integrals that are transformers of $\bS_p$ into itself.

The notion of Schur multipliers can be generalized to the case of spectral measures on Hilbert space (this terminology was introduced in \cite{Pe}). To define Schur multipliers with respect to spectral measures, we are going to use double operator integrals.

{\it Double operator integrals} are expressions of the form
$$
\iint\limits_{\X_1\times\X_2}\Phi(x,y)\,dE_1(x)Q\,dE_2(y),
$$
where $E_1$ and $E_2$ are spectral measures, $Q$ is a liner operator no Hilbert space and $\Phi$ is a bounded measurable function.

Double operator integrals appeared first in the paper \cite{DK}. Later Birman and Solomyak elaborated a beautiful theory of double operator integrals and found important applications, see \cite{BS1}. We refer the reader to recent surveys \cite{Pe2} and \cite{AP2} for definitions and basic properties of double operator integrals.

\medskip

{\bf Definition.}
The function $\Phi$ is called a {\it Schur multiplier of $\bS_p$ with respect to the spectral measures
$E_1$ and $E_2$} if
$$
\iint\Phi(x,y)\,dE_1(x)Q\,dE_2(y)\in\bS_p\quad\mbox{whenever}\quad Q\in\bS_p.
$$
We use the notation $\fM_p({E_1,E_2})$ for the class of Schur multipliers with respect to $E_1$ and $E_2$.

\medskip

Let us explain that this is a generalization of the notion of a matrix Schur multiplier. Indeed,
let $E_1$ and $E_2$ be the spectral measures on $\ell^2$ defined by
\bay
\label{smonl2}
E_1(\D)v=E_2(\D)v\df\sum_{j\in\D}(v,\ee_j)\ee_j,\quad\D\subset\Z_+,\quad v\in\ell^2,
\ey
where $\{\ee_j\}_{j\ge0}$ is the standard orthonormal basis of $\ell^2$.
It is easy to see that if $Q$ is a linear operator on $\ell^2$ with matrix 
$\{q_{jk}\}_{j,k\ge0}$ and $\Phi=\{\Phi(j,k)\}_{j,k\ge0}$ is a matrix with bounded entries, then the double operator integral
$$
\iint\limits_{\Z_+\times\Z_+}\Phi(j,k)\,dE_1(j)Q\,dE_2(k)
$$
is a linear operator with matrix $\{\Phi(j,k)q_{jk}\}_{j,k\ge0}$. Therefore, 
$\Phi\in\fM_p({E_1,E_2})$ if and only if the matrix 
$\{\Phi(j,k)\}_{j,k\ge0}\in\fM_p$.

To state a well-known sufficient condition for the membership in the space of Schur multipliers of $\bS_p$, $p\le1$, we introduce the notion of $p$-tensor product of $L^\be$ spaces.

\medskip

{\bf Definition.} Let $0<p\le1$ and let $E_1$ and $E_2$ be spectral measures defined on $\s$-algebras of subsets of $\X_1$ and $\X_2$. The $p$-tensor product
$L^\be_{E_1}\otimes_pL^\be_{E_2}$ is, by definition, the class of functions $\Phi$ on
$\X_1\times\X_2$ of the form 
\bay
\label{tensu}
\Phi(x,y)=\sum_{n\ge0}\f_n(x)\psi_n(y),
\ey
where the functions $\f_n$ in $L^\be_{E_1}$ and $\psi_n$ in $L^\be_{E_2}$ satisfy
\bay
\label{pner}
\left(\sum_{n\ge0}\|\f_n\|^p_{L^\be}\|\psi_n\|^p_{L^\be}\right)^{1/p}<\be.
\ey
By $\|\Phi\|_{L^\be\otimes_p L^\be}$ we mean the infimum of the left-hand side of 
\rf{pner} over all representations of $\Phi$ in the form \rf{tensu}.

\medskip 

Note that in the case $p=1$, the space $L^\be_{E_1}\otimes_1L^\be_{E_2}$
is the projective tensor product $L^\be_{E_1}\widehat\otimes L^\be_{E_2}$ of
$L^\be_{E_1}$ and $L^\be_{E_2}$.

It is well known and it is easy to see that 
$L^\be_{E_1}\otimes_pL^\be_{E_2}\subset\fM_p({E_1,E_2})$
and 
\bay
\label{odnomen'shedrugogo}
\|\Phi\|_{\fM_p({E_1,E_2})}\le\|\Phi\|_{L^\be\otimes_p L^\be},
\quad\Phi\in L^\be_{E_1}\otimes_pL^\be_{E_2}.
\ey
Indeed, suppose that
$\Phi$ is given by \rf{tensu} and \rf{pner} holds and let $Q\in\bS_p$. We have 
\begin{align*}
\left\|\iint\Phi\,dE_1Q\,dE_2\right\|^p_{\bS_p}&=
\left\|\sum_{n\ge0}\left(\int\f_n\,dE_1\right)Q\left(\int\psi_n\,dE_2\right)
\right\|^p_{\bS_p}\\[.2cm]
&\le\sum_{n\ge0}\|\f_n\|_{L^\be}^p\|\psi_n\|_{L^\be}^p\|Q\|^p_{\bS_p}.
\end{align*}

Let us consider the space $L^\be_{E_1}\otimes^{\rm c}_pL^\be_{E_2}$, which consists of functions $\Phi$ on $\X_1\times\X_2$, for which there exists a sequence of functions 
$\Phi_n$ in $L^\be_{E_1}\otimes_pL^\be_{E_2}$ such that
$$
\|\Phi_n\|_{L^\be_{E_1}\otimes_pL^\be_{E_2}}\le\const\quad\mbox{and}\quad
\Phi_n(x,y)\to\Phi(x,y)\quad\mbox{a.e.}
$$

In particular, we can consider the $p$-tensor product $\ell^\be\otimes_p\ell^\be$, which corresponds to the case when $E_1$ and $E_2$ are the spectral measures on $\ell^2$ defined by \rf{smonl2}. It is clear from the definition of $(r,p,q)$-nuclear operators (see \cite{Pie}, \S18.1.1) that a matrix belongs to  
$\ell^\be\otimes_p\ell^\be$, $0<p\le1$, if and only if it induces a $(p,1,1)$-nuclear operator
from $\ell^1$ to $\ell^\be$. Note that $(p,1,1)$-nuclear operators and $p$-projective tensor products  for $p\in(0,1]$ were
investigated by Grothendieck \cite{Gro},
Chapter 2.

In the case when $E_1$ and $E_2$ are the spectral measures on $\ell^2$ defined by \rf{smonl2}, the space $\ell^\be\otimes_p^{\rm c}\ell^\be$ is the space of infinite matrices
$A$ such that 
$$
\|P_nAP_n\|_{\ell^\be\otimes_p\ell^\be}\le\const,
$$
where $P_n$ is the standard natural projection on $\ell^\be$ onto 
$\spn\{\ee_j:~0\le j\le n\}$.

It is easy to see that $L^\be_{E_1}\otimes^{\rm c}_pL^\be_{E_2}$ is still a subset of 
$\fM_p({E_1,E_2})$, $0<p\le1$. In particular,
\bay
\label{perdous}
\ell^\be\otimes_p^{\rm c}\ell^\be\subset\fM_p,\quad0<p\le1.
\ey
It is well known that $L^\be_{E_1}\otimes^{\rm c}_1L^\be_{E_2}=\fM_1({E_1,E_2})$, see
\cite{Pe}. 

It is also well known that 
$L^\be_{E_1}\otimes_1L^\be_{E_2}=L^\be_{E_1}\widehat\otimes L^\be_{E_2}$ in general is a proper subset of $\fM_1({E_1,E_2})$. Indeed, consider the case when
$E_1$ and $E_2$ are the spectral measures on $\ell^2$ defined by \rf{smonl2}. Consider the identity matrix $I$, i.e., the infinite matrix with diagonal entries equal to 1 and off diagonal entries equal to 0. Obviously, $I\in\fM_1=\ell^\be\otimes_1^{\rm c}\ell^\be$. On the other hand,  $I\not\in\ell^\be\widehat\otimes\ell^\be$. Indeed, the space 
$\ell^\be\widehat\otimes\ell^\be$ can naturally be identified with the nuclear operators
(we refer the reader to \cite{Pie} for the definition of nuclear operators) from $\ell^1$ to $\ell^\be$, and so if $I$ belonged to $\ell^\be\widehat\otimes\ell^\be$, it would induce a compact operator from $\ell^1$ to $\ell^\be$. However, it is not compact. We can consider the standard basis $\{\ee_j\}_{j\ge0}$ in $\ell^1$. Clearly, there is no subsequence of
$\{\ee_j\}_{j\ge0}$ in $\ell^1$ that would converge in the norm of $\ell^\be$.

However, {\it we do not know} whether for $p<1$, 
$$
L^\be_{E_1}\otimes_pL^\be_{E_2}=L^\be_{E_1}\otimes^{\rm c}_pL^\be_{E_2}.
$$ 

It will be shown in \S\:\ref{pqr} that for $p<1$, 
$$
\ell^\be\otimes^{\rm c}_p\ell^\be\ne\fM_p.
$$
This is a combination of an unpublished observation by Nigel Kalton and the description of the class of diagonal matrices in $\fM_p$ obtained in \cite{AP}.

In \S\:\ref{uzhas} we introduce a scale $\W_p$ of tensor products of two $\ell^\be$ spaces and prove that the matrices in $\W_p$ are Schur multipliers of $\bS_p$ for $p<1$. We compare this sufficient condition with the sufficient condition in terms of 
the space $\ell^\be\otimes^{\rm c}_p\ell^\be$.

We analyze in \S\:\ref{Pisier} another sufficient condition for a matrix to be a Schur multiplier of $\bS_p$ for $p<1$. This sufficient condition follows from a result of Pisier \cite{Pis2}. We compare this condition with the other sufficient conditions discussed in this paper.

Let us also mention that in \cite{AP} we studied two special classes of matrix Schur multipliers: Toeplitz Schur multipliers and Hankel Schur multipliers.

Throughout this paper we consider only separable Hilbert spaces.

\

\section{\bf Schur multipliers of $\bs{\bS_p}$ with respect to spectral measures\\ and integral operators}
\setcounter{equation}{0}
\label{dois}

\

Suppose that $(\X_1,{\frak B}_1)$ and $(\X_2,{\frak B}_2)$ are measurable spaces and  $E_1$ and $E_2$ are spectral measure on ${\frak B}_1$ and ${\frak B}_2$ that take values in the set of orthogonal projections on a Hilbert space $\h$. Let $\mu_1$ be a positive measure that is mutually absolutely continuous with $E_1$ and let 
$\mu_2$ be a positive measure that is mutually absolutely continuous with $E_2$.

For a function $k$ in $L^2(\mu_1\otimes\mu_2)$, we consider the integral operator 
${\mathcal I}_k:L^2(\mu_2)\to L^2(\mu_1)$ defined by
\bay
\label{Ik}
\big({\mathcal I}_kf\big)(x)=\int_{\X_2}k(x,y)f(y)\,d\mu_2(y),\quad f\in L^2(\mu_2).
\ey

\begin{thm}
\label{Schur=Schur}
Let $0<p\le1$ and let
$\Phi$ be a bounded measurable functions on $\X_1\times\X_2$. The following are equivalent:

{\em(i)} $\Phi\in\fM_p({E_1,E_2})$;

{\em(ii)} if $k$ is a function on $\X_1\times\X_2$ in $L^2(\mu_1\otimes\mu_2)$
such that ${\mathcal I}_k\in\bS_p$, then ${\mathcal I}_{\Phi k}\in\bS_p$.
\end{thm}

In other words, statement (ii) of Theorem \ref{Schur=Schur} means that $\f$ {\it is a multiplier of the space of kernel functions of integral operators of class} $\bS_p$.

Note that the case $p=1$ is well-known, see \cite{BS1}, \cite{Pe}.

We need an auxiliary fact known to experts. Let $E$ be a spectral measure on a Hilbert space $\h$ and let $u$ and $v$ be vectors in $\h$. Consider the positive measure $\s_u$ and the complex measure $\s_{u,v}$ defined by $\s_u(\D)=(E(\D)u,u)$ and the complex measure $\s_{u,v}$
defined by
$$
\s_u(\D)=(E(\D)u,u),\quad \s_{u,v}(\D)=(E(\D)u,v),
$$
where $\D$ is a measurable set. Since 
$$
|\s_{u,v}(\D)|\le\|E(\D)u\|\cdot\|v\|=(E(\D)u,u)^{1/2}\|v\|,
$$
it follows that $\s_{u,v}$ is absolutely continuous with respect to $\s_u$. We denote by $h_v$ the Radon--Nikodym density of $\s_{u,v}$, i.e., $d\s_{u,v}=h_vd\s_u$.

\begin{lem}
\label{hvsu}
For an arbitrary vector $v$ in $\h$, the function $h_v$ belongs to $L^2(\s_u)$.
Moreover,
$$
\big\{h_v:~v\in\h\big\}=L^2(\s_u).
$$
\end{lem}

Lemma \ref{hvsu} can easily be deduced from the theory of spectral multiplicities, see \cite{BS3}, Ch. 7. However, for the sake of convenience, we give a more elementary proof here.

\medskip

\Pf Let $h\in L^2(\s_u)$. Let $A$ be the (not necessarily bounded) operator defined by
$
A=\int\ov h\,dE
$
whose domain ${\rm D}(A)$ is equal to
$$
{\rm D}(A)=\left\{w\in\h:~\int|h|^2\,d\s_w<\be\right\},
$$
and so $u\in{\rm D}(A)$. Put $v\df Au$. It is easy to see that 
$$
\big(E(\D)u,v\big)=\int_\D h\,d\s_u,
$$
and so $h_v=h$.

Suppose now that $v\in\h$. Let us show that $h_v\in L^2(\s_u)$. Let $g$ be a simple function, i.e., a linear combination of characteristic functions. 
Then
$$
\int h_vg\,d\s_u=\left(\left(\int g\,dE\right)u,v\right)
$$
(it suffices to verify this identity for characteristic functions $g$). Hence,
$$
\left|\int h_vg\,d\s_u\right|\le\|v\|\cdot\left\|\left(\int g\,dE\right)u\right\|
\le\|v\|\cdot\|g\|_{L^2(\s_u)}.
$$
Since this inequality holds for an arbitrary simple function $g$, it follows that 
$h_v\in L^2(\s_u)$. $\bl$

\medskip

{\bf Remark.}
Before we proceed to the proof of Theorem \ref{Schur=Schur}, let us make the following observation. Suppose that $\mu_1^\flat$ is a positive measure that is mutually absolutely continuous with $\mu_1$ and $\mu_2^\flat$ is a positive measure that is mutually absolutely continuous with $\mu_2$. Then the class of multipliers of the space of kernel functions of integral operators from $L^2(\mu_2^\flat)$ to $L^2(\mu_1^\flat)$ of class $\bS_p$ coincides with the class of multipliers the space of kernel functions of integral operators from $L^2(\mu_2)$ to $L^2(\mu_1)$ of class $\bS_p$. This can be seen very easily by considering unitary operators 
$U:L^2(\mu_1^\flat)\to L^2(\mu_1)$ and  $V:L^2(\mu_2^\flat)\to L^2(\mu_2)$
defined by
$$
Uf=h_1^{1/2}f,\quad f\in L^2(\mu_1^\flat),\qquad\mbox{and}\qquad 
Vf=h_2^{1/2}f,\quad f\in L^2(\mu_2^\flat),
$$
where $d\mu_1^\flat=h_1\,d\mu_1$ and $d\mu_2^\flat=h_2\,d\mu_2$.

\medskip

{\bf Proof of Theorem \ref{Schur=Schur}.} Put
$$
Q_{u,v}=\iint_{\X_1\times\X_2}\Phi(x,y)\,dE_1(x)T_{u,v}\,dE_2(y),
$$
where $T$ is the rank one operator $T=(\cdot,u)v$, $u,\,v\in\h$. Since $\Phi$ is a bounded function, $Q\in\bS_2$.

Clearly, $\Phi$ is a Schur multiplier with respect to $E_1$ and $E_2$ is and only if 
$Q_{u,v}\in\bS_p$ and $\|Q_{u,v}\|_{\bS_p}\le\const\|u\|\cdot\|v\|$ for arbitrary vectors $u$ and $v$.

It is easy to verify (see \cite{BS1}) that
$$
(Q_{u,v}w,z)=\iint \Phi(x,y)\,d\s_{w,u}(x)\,d\l_{v,z}(y),
$$
where
$$
\s_{w,u}(\D_2)\df\big(E_2(\D_2)w,u\big)\quad\mbox{and}\quad
\l_{v,z}(\D_1)\df\big(E_1(\D_1)v,z\big),
$$
for arbitrary measurable subsets $\D_1$ and $\D_2$ of $\X_1$ and $\X_2$. 
Therefore, taking into account Lemma \ref{hvsu}, 
we can conclude that statement (i) is equivalent to the fact that
the integral operator 
${\mathcal I}_{\Phi(h\otimes g)}$ with kernel function
$(x,y)\mapsto\Phi(x,y)h(x)g(y)$  belongs to $\bS_p$ for arbitrary $g$ in $L^2(\s_v)$ and $h$ in $L^2(\s_u)$ and
$$
\big\|{\mathcal I}_{\Phi(h\otimes g)}\big\|_{\bS_p}\le\const
\|g\|_{L^2(\s_v)}\|h\|_{L^2(\s_u)}
$$
for arbitrary $g$ in $L^2(\s_v)$ and $h$ in $L^2(\s_u)$. This exactly means that
the function $\Phi$ is a multiplier of the space kernel functions of $S_p$ integral operators from $L^2(\s_v)$ to $L^2(\s_u)$.

To complete the proof, it remains to find vectors $u$ and $v$ such that $\s_u$ is mutually absolutely continuous with $E_2$ and $\s_v$ is mutually absolutely continuous with $E_1$. The existence of such vectors is well known, see e.g., \cite{BS3}, Ch. 7. $\bl$

\

\section{\bf Completely bounded matrix Schur multipliers}
\setcounter{equation}{0}
\label{vpoogr}

\

We prove in this section that for $p\le1$, matrix Schur multipliers of $\bS_p$ are necessarily completely bounded.

\begin{thm}
Let $0<p\le1$. Suppose that a matrix $A=\{\a_{jk}\}_{j,k\ge0}$ is a Schur multiplier of $\bS_p$. Then $A$ is a completely bounded Schur multiplier of $\bS_p$.
\end{thm}

\Pf Let $\K$ be a Hilbert space. Consider the spectral measure $E$ on $\ell^2(\K)$
defined by 
$$
E(\D)\{u_n\}_{n\ge0}=\{\Bbbone_\D(n)u_n\}_{n\ge0}
$$
where $\D$ is a subset of $\Z_+$ and $\Bbbone_\D$ is its characteristic function.
Let $\nu$ be the counting measure on $\Z_+$, i.e., $\nu(\D)$ is the number of elements of $\D$. Clearly, $\nu$ is mutually absolutely continuous with $E$.

It is easy to see the condition that a matrix $A=\{\a_{jk}\}_{j,k\ge0}$ is 
a matrix Schur multiplier of $\bS_p$ exactly means that it is a multiplier of the 
space of kernel functions of $\bS_p$ integral operators on $L^2(\nu)$.

On the other hand, it is easy to see that if $T$ is an operator on $\ell^2(\K)$ with block matrix $\{T_{jk}\}_{j,k\ge0}$, then
$$
\iint\limits_{\Z_+\times\Z_+}\a_{jk}\,dE(j)T\,dE(k)
$$
is the operator on $\ell^2(\K)$ with block matrix $\{\a_{jk}T_{jk}\}_{j,k\ge0}$.

The result follow now from Theorem \ref{Schur=Schur}. $\bl$

\medskip

To conclude this section we state the following problem.

\medskip

{\bf Problem.} Can Theorem \ref{Schur=Schur} be generalized to the case $p>1$?

\medskip

Note that an affirmative answer to this question would imply that matrix Schur multipliers of $\bS_p$ with $p>1$ must be completely bounded.

\

\section{\bf Diagonal matrices in $\bs{\ell^\be\otimes_p\ell^\be}$}
\setcounter{equation}{0}
\label{pqr}

\

In \S\:\ref{In} for $p\le1$, we defined the $p$-tensor product $\tp$. In this section we are going to study diagonal matrices in $\tp$. Recall that $\tp$ consists of infinite matrices of the form
$$
\sum_{n\ge0} x_n\otimes y_n\df\sum_{n\ge0} x_n y_n^{\rm t}
$$
where $x_n$ and $y_n$ are columns in $\ell^\be$ such that
$$
\sum_{n\ge0}\|x_n\|^p_{\ell^\be}\|y_n\|_{\ell^\be}^p<\be.
$$
Here and in what follows, for a matrix $C$, we use the notation $C^{\rm t}$ for the transposed matrix.

%
%
%
%

Following \cite{AP}, we associate  with $p\in(0,1]$ the number $p_\sharp$, $p_\sharp\in(0,\be]$, defined by
$$
\ps\df\frac{p}{1-p}.
$$

The following theorem is a special case of Theorem 3.2 in \cite{AP}.

\begin{thm}
\label{diag}
Let $p\le1$ and let $M$ be an infinite diagonal matrix with diagonal entries $\mu_j$, $j\ge0$. Then
$M\in\fM_p$ if and only if  $\{\mu_j\}_{j\ge0}\in\ell^{\ps}$. Moreover,
$\|M\|_{\fM_p}=\|\mu\|_{\ell^{\ps}}=\|M\|_{\bS_{\ps}}$.
\end{thm}

{\bf Definition.}
Let $p\in(0,1]$. Denote by $D_p$ the space of sequences $\mu=\{\mu_j\}_{j\ge0}$
such that the diagonal matrix $M$ with diagonal entries $\{\mu_j\}_{j\ge0}$
belongs to $\tp$.
Put $\|\mu\|_{D_p}\df\|M\|_{\ell^\be\otimes_p\ell^\be}$.

\medskip

Since $\tp\subset\fM_p$, it follows from Theorem  \ref{diag} that $D_p\subset\ell^{\ps}$.

Let $\mu=\{\mu_j\}_{j\ge0}$ be a sequence whose terms tend to 0.
Denote by $\mu^*=\{\mu^*_j\}_{j\ge0}$ the nonincreasing rearrangement of
the sequence $\{|\mu_j|\}_{j\ge0}$.

Let us remind the definition  of the Lorentz space $\ell^{q,r}$ for $r,q\in(0,\be]$:
$$
\ell^{q,r}\df\left\{\mu\in c_0: \|\mu\|^r_{\ell^{q,r}}\df\sum_{j\ge0}(\mu_j^*)^r(1+j)^{\frac rq-1}<\be\right\},\quad r<\be,
$$
$$
\ell^{q,\be}\df\left\{\mu\in c_0: \|\mu\|_{\ell^{q,\be}}\df\sup_{j\ge0}(\mu_j^*)(1+j)^{\frac 1q}<\be\right\}\quad\mbox{and}\quad\ell^{\be,\be}\df\ell^\be.
$$

Let $\bS_{q,r}$ denote the corresponding operator ideal in the space of bounded operators on Hilbert space, i.e.,
$$
\bS_{q,r}=\big\{T:~\{s_j(T)\}_{j\ge0}\in\ell^{q,r}\big\},
$$
where $\{s_j(T)\}_{j\ge0}$ is the sequence of singular values of $T$.

\begin{thm}
\label{diag2}
Let $0<p\le1$. Then $\ell^{p_\sharp,p}\subset D_p$, i.e., if 
$\sum\limits_{j\ge0}\frac{(\mu_j^*)^p}{(1+j)^p}<\be$,
then $\{\mu_j\}_{j\ge0}\in D_p$. Moreover,
$$
\|\mu\|_{D_p}^p\le2\sum_{j\ge0}\frac{(\mu_j^*)^p}{(1+j)^p}<\be.
$$
\end{thm}


\begin{lem}
\label{edindiag}
Let $p\in(0,1]$ and
let $J_n=\sum_{k=0}^{n-1}\ee_k\otimes\ee_k$. Then $\|J_n\|_{\tp}=n^{1/p_\sharp}$.
\end{lem}

\Pf
Put 
$$
x_k=\sum_{j=0}^{n-1}e^{\frac{2\pi{\rm i}kj}n}\ee_j\quad\mbox{and}\quad 
y_k=\sum_{j=0}^{n-1}e^{-\frac{2\pi{\rm i}kj}n}\ee_j,\qquad0\le k\le n-1.
$$
It can be easily verified that $J_n=n^{-1}\sum_{k=0}^{n-1}x_k\otimes y_k$. Hence,
$\|J_n\|_{\ell^\be\otimes_p\ell^\be}\le n^{1/\ps}$.
It remains to observe that $\|J_n\|_{\fM_p}\le\|J_n\|_{\ell^\be\otimes_p\ell^\be}$
by \rf{odnomen'shedrugogo}  and $\|J_n\|_{\fM_p}=n^{1/\ps}$
by Theorem \ref{diag}. $\bl$

\begin{lem} 
\label{razshur1}
Let $p\in(0,1]$. For matrices $M=\{\mu_{jk}\}_{j,k\ge0}$ and $A=\{a_{jk}\}_{j,k\ge0}$, the following inequality holds:
$$
\|AM\|_{\ell^\be\otimes_p\ell^\be}\le\|A\|_{\ell^\be\to\ell^\be}\|M\|_{\ell^\be\otimes_p\ell^\be},
$$
where $\|A\|_{\ell^\be\to\ell^\be}$ denote the operator norm on $\ell^\be$.
\end{lem}

\Pf Let $M\in\ell^\be\otimes_p\ell^\be$ and let $\e>0$. Then $M$ can be represented in the form
$$
M=\sum_{n\ge0} x_n y_n^{\rm t},
$$
where $x_n$ and $y_n$ are columns in $\ell^\be$ such that
$$
\sum_{n\ge0}\|x_n\|^p_{\ell^\be}\|y_n\|_{\ell^\be}^p<\|M\|^p_{\ell^\be\otimes_p\ell^\be}+\e.
$$
It follows that
$$
AM=\sum_{n\ge0}Ax_n y_n^{\rm t}
$$
and
$$
\|AM\|_{\ell^\be\otimes_p\ell^\be}^p\le\sum_{n\ge0}\|Ax_n\|^p_{\ell^\be}\|y_n\|_{\ell^\be}^p
\le \|A\|_{\ell^\be\to\ell^\be}^p(\|M\|_{\ell^\be\otimes_p\ell^\be}^p+\e). \quad \bl
$$

\begin{cor} 
\label{razdiag1}
Let $0<p\le1$ and
let $A$ and $M$ be infinite diagonal matrices with diagonal entries
$\{\a_j\}_{j\ge0}$ and $\{\mu_j\}_{j\ge0}$. Then 
$$
\|AM\|_{\ell^\be\otimes_p\ell^\be}\le\Big(\sup_{j\ge0}|\a_j|\Big)\|M\|_{\ell^\be\otimes_p\ell^\be}.
$$
\end{cor}

{\bf Remark.} In the same way one can prove that
$$
\|MA\|_{\ell^\be\otimes_p\ell^\be}\le\|A^{\rm t}\|_{\ell^\be\to\ell^\be}\|M\|_{\ell^\be\otimes_p\ell^\be}=
\|A\|_{\ell^\be\to\ell^\be}\|M\|_{\ell^1\otimes_p\ell^1},
$$
where $\|A\|_{\ell^1\to\ell^1}$ denote the operator norms on $\ell^1$.

\medskip

Recall that
$$
\|A\|_{\ell^1\to\ell^1}=\sup\limits_k\sum\limits_{j\ge0}|a_{jk}|\quad
\mbox{and}\quad 
\|A\|_{\ell^\be\to\ell^\be}=\sup\limits_j\sum\limits_{k\ge0}|a_{jk}|.
$$

\medskip

{\bf Proof of Theorem \ref{diag2}.} By Corollary \ref{razdiag1}, it suffices to consider the case when
$\mu_j\ge0$ for all $j\ge0$. Moreover, we may assume that $\{\mu_j\}_{\ge0}$ is
a decreasing sequence. Put $M_n=M\star J_n$, where $J_n$ denotes the same as in
Lemma \ref{edindiag}.
Lemma \ref{edindiag} and Corollary \ref{razdiag1} imply that 
$\|M_{2^{n+1}}-M_{2^n}\|_{\ell^\be\to\ell^\be}\le2^{n/\ps}\mu_{2^n}$.
It follows that
\bey
\|M\|_{\ell^\be\otimes_p\ell^\be}^p\le\|M_1\|_{\ell^\be\otimes_p\ell^\be}^p+
\sum_{n\ge0}\|M_{2^{n+1}}-M_{2^n}\|_{\ell^\be\otimes_p\ell^\be}^p
\le\mu_0^p+\sum_{n\ge0}2^{pn/\ps}\mu_{2^{n}}^p\\
=\mu_0^p+\sum_{n\ge0}2^n\frac{\mu_{2^n}^p}{2^{pn}}
\le\mu_0^p+\mu_1^p+2\sum_{j\ge1}\frac{\mu_j^p}{(1+j)^p}
\le2\sum_{j\ge0}\frac{\mu_j^p}{(1+j)^p}.\quad \bl
\eey

\

\section{$\bs{\ell^\be\otimes^{\rm c}_p\ell^\be\ne\fM_p}$}
\setcounter{equation}{0}
\label{diagon}

\

The main purpose of this section is to show that for $p<1$, the space 
$\ell^\be\otimes^{\rm c}_p\ell^\be$ is a proper subset of $\fM_p$.

We can naturally imbed the finite-dimensional space $\C^n$ to $\ell^\be$ and we are going to use the notation  $\|x\|_{\ell^\be}$ for vectors $x$ in $\C^n$ and 
$\|A\|_{\ell^\be\otimes_p\ell^\be}$ for $n\times n$ matrices $A$.
As usual, $I_n$ stands for the $n\times n$ identity matrix.

\begin{lem} 
\label{l1}
Let $p\in(0,1)$. Suppose that
$$
n^{-1/{p_\sharp}}I_n = \sum_{k\ge1}\l_k u_k\otimes v_k,
$$
where $u_k,v_k\in\C^n$ with $\|u_k\|_{\ell^\be},\|v_k\|_{\ell^\be}\le 1$,
$\{\l_k\}_{k\ge1}$ is a nonincreasing summable sequence of
nonnegative numbers. Let $a$ and $b$ be positive numbers such that $a<1<b$. Then
$$
\sum_{an<k\le bn} \l_k^p \ge (1-a)
-b^{p-1}\|\l\|_{\ell_p}^p.
$$
\end{lem}

\Pf
Note that the rank of the operator $\sum\limits_{k\le an}\l_k u_k\otimes v_k$ is at most $an.$  Hence,
the operator $n^{1/{p_\sharp}}\sum\limits_{k>an}\l_k u_k\otimes v_k$ is identical on a subspace of dimension at least $(1-a)n$, and so 
$$
n^{1-p}\left\|\sum_{k>an}\l_ku_k\otimes v_k\right\|_{\bS_p}^p
 \ge n(1-a).
$$
Moreover, 
\begin{align*} 
\left\|\sum_{k>bn} \l_ku_k\otimes v_k\right\|_{\bS_p}&\le
n^{1/p_\sharp}\left\|\sum_{k>bn}\l_k u_k\otimes v_k\right\|_{\bS_1}
\le n^{1/p_\sharp}\sum_{k>bn} |\l_k |\cdot\|u_k\otimes v_k\|_{\bS_1}\\[.2cm]
&\le n^{\frac1p}\sum_{k>bn} |\l_k |
\le n^{1/p}\Big(\max_{k>bn}\l_k \Big)^{1-p}\sum_{k>bn} \l_k^p
\le b^{-1/p_\sharp}n\|\l\|_{\ell^p}.
\end{align*}
Hence,
\begin{align*} 
n^p\sum_{an<k\le bn}\l_k^p&\ge
\left\|\sum_{an<k\le bn}\l_ku_k\otimes v_k\right\|_{\bS_p}^p\\[.2cm]
& \ge
 \left\|\sum_{k>an}\l_ku_k\otimes v_k\right\|_{\bS_p}^p
 -\left\|\sum_{k>bn}\l_ku_k\otimes v_k\right\|_{\bS_p}^p
 \\[.2cm]
 &\ge n^p(1-a)
-b^{p-1}n^p\|\l\|_{\ell_p}^p.\quad\bl
\end{align*}

Let $\{A_j\}_{j=1}^k$ be a sequence of matrices of size $n_j\times n_j$.  We can naturally identify  $\bigoplus_{j=1}^k A_j$ with a matrix of size
$\big(\sum_{j=1}^kn_j\big)\times\big(\sum_{j=1}^kn_j\big)$.

\begin{lem}
\label{l2}
Let $p\in(0,1)$  and
let $\{n_j\}_{j=1}^m$ be a finite sequence of positive integers such that
$n_{j+1}\ge4(2m)^{\frac1{1-p}}n_j$. Then
$$
\left\|\bigoplus_{j=1}^mn_j^{-1/p_\sharp}I_{n_j}\right\|^p_{\ell^\be\otimes_p\ell^\be}\ge \frac m2.
$$
\end{lem}

\Pf Let
$$
\bigoplus_{j=1}^mn_j^{-1/p_\sharp}I_{n_j}= \sum_{k=1}^\be\l_k u_k\otimes v_k,
$$
where $\|u_k\|_{\be},\|v_k\|_{\be}\le 1$ and
$\{\l_k\}_{k\ge1}$ is a nonincreasing summable sequence of
nonnegative numbers. We have to show that $\|\l\|_{\ell^p}^p\ge \frac m2$.
Applying Lemma \ref{l1} for $a=\frac14$ and $b=(2m)^{\frac1{1-p}}$, we get
$$
\sum_{an_j<k\le bn_j} \lambda_k^p \ge\frac34-\frac1{2m}\|\l\|_{\ell_p}^p
$$
whenever $1\le j\le m$.
Keeping in mind that $bn_j\le an_{j+1}$, we obtain
$$
\|\l\|_{\ell_p}^p\ge
\sum_{k=1}^\be \lambda_k^p \ge\sum_{j=1}^m\left(\sum_{an_j<k\le bn_j} \lambda_k^p\right)
\ge\frac34m-\frac12\|\l\|_{\ell_p}^p.
$$
Hence, $\|\l\|_{\ell^p}^p\ge \frac m2$. $\bl$

\begin{thm} 
\label{16}
Let $0<p<1$. Then $\fM_p\not=\ell^\be\otimes_p\ell^\be$. 
Moreover, there exists
a diagonal matrix $M$ such that $M\in\fM_p$ but $M\notin\ell^\be\otimes_p\ell^\be$.
In other words, $\ell^{p_\sharp}\ne D_p$.
\end{thm}

\Pf Assume that each diagonal matrix $M$ in $\fM_p$ belongs to the space $\ell^\be\otimes_p\ell^\be$. Then by the closed graph
theorem, there exists a constant $c$ such that 
$\|M\|_{\ell^\be\otimes_p\ell^\be}\le c\|M\|_{\fM_p}$. By Lemma \ref{l2}, we have
\bay
\label{niz}
\left\|\bigoplus_{j=1}^mn_j^{-1/p_\sharp}I_{n_j}\right\|_{\ell^\be\otimes_p\ell^\be}\ge
2^{-1/p}m^{1/p}.
\ey
On the other hand, by Theorem \ref{diag}, 
\bay
\label{verh}
\left\|\bigoplus_{j=1}^mn_j^{-1/p_\sharp}I_{n_j}\right\|_{\fM_p}
=\left\|\bigoplus_{j=1}^mn_j^{-1/p_\sharp}I_{n_j}\right\|_{\bS_{p_\sharp}}=m^{1/p_\sharp}
\ey
and we get a contradiction. $\bl$

\medskip

{\bf Remark.}
It is easy to see that in the notation of the proof of Theorem \ref{16} we have
$$
\left\|\bigoplus_{j=1}^mn_j^{-1/p_\sharp}I_{n_j}\right\|_{\bS_{p_\sharp,r}}\le\const m^{1/r}.
$$

\medskip

It is clear from this remark that the proof of  Theorem \ref{16} allows us to get the following
fact, which complements Theorem \ref{diag2}.


\begin{thm}
\label{diag2kontra}
Let $0<p<1$ and let $r>p$. Then $\ell^{p_\sharp,r}\not\subset D_p$ .
%
\end{thm}

\medskip

{\bf Remark.} As we have mentioned in \S\:\ref{In}, the identity matrix does not belong to the space $\ell^\be\otimes_1\ell^\be$.
Hence, Theorem  \ref{16} remains true for $p=1$. On the other hand, Theorem \ref{diag2kontra} cannot be generalized to the case 
$p=1$ and $r<\be$ because in this case $\ell^{p_\sharp,r}\subset c_0=D_1$.

\medskip

%


\medskip

For $p\in(0,1]$, we denote by $D_p^{\rm c}$ the space of all sequences $\mu=\{\mu_j\}_{j\ge0}$
such that the diagonal matrix $M$ with diagonal entries $\{\mu_j\}_{j\ge0}$
belongs to $\ell^\be\otimes^{\rm c}_p\ell^\be$.
Clearly, $D^{\rm c}_p\subset\ell^{\ps}$.

Theorem \ref{16} allows us to prove the following fact.

\begin{thm} 
\label{Dpll}
Let $p\in(0,1)$. Then $D_p^{\rm c}\ne\ell^{\ps}$.
\end{thm}

\Pf Suppose that $D_p^{\rm c}=\ell^{\ps}$. Then there exists a constant $C(p)$ such that 
$\|\mu\|_{D_p^{\rm c}}\le C(p)\|\mu\|_{\ell^{\ps}}$ 
for all sequences $\mu=\{\mu_j\}_{j\ge0}$. Note that $\|\mu\|_{D_p^{\rm c}}=\|\mu\|_{D_p}$
if $\mu_j=0$ for all sufficient large $j$. Thus, we get a contradiction with inequalities
\rf{verh} and \rf{niz}. $\bl$

\medskip

Theorem \ref{Dpll} immediately implies the main result of this section.

\begin{thm}
$\ell^\be\otimes^{\rm c}_p\ell^\be\ne\fM_p$.
\end{thm}

To conclude the section, we state the following conjecture:

\medskip

{\bf Conjecture.} Let $p<1$. Suppose that $M$ is a diagonal matrix. Then the following are equivalent:

(i) $M\in\ell^\be\otimes_p\ell^\be$;

(ii) $M\in\ell^\be\otimes^{\rm c}_p\ell^\be$;

(iii) $\|M\|_{\bS_{\ps,p}}<\be$.

\

\section{\bf Another scale of tensor products}
\setcounter{equation}{0}
\label{uzhas}

\

In this section we introduce a scale $\W_p$,  $0<p\le1$, of tensor products of $\ell^\be$ spaces. The right endpoint of this scale coincides with the Haagerup tensor product of $\ell^\be$ spaces. We show that  $\W_p\subset\fM_p$. 

This together with \rf{perdous} gives us two sufficient conditions for a matrix to belong to $\fM_p$. We show in this section that for $p<1$, none of them implies the other one.

For $p\in(0,1)$, we put $\Xi_p\df\ell^{2p_\sharp}(\ell^{2p})$. In other words, $\Xi_p$ is the set of matrices 
$X=\{x_{jk}\}_{j,k\ge0}$ such that 
$$
\|X\|_{\Xi_p}\df\left(\sum_{j\ge0}\left(\sum_{k\ge0}|x_{jk}|^{2p}\right)^{\ps/p}\right)^{\frac1{2\ps}}<\be.
$$

For $p=1$, put $\Xi_1\df\ell^{\be}(\ell^{2})$ and
$$
\|X\|_{\Xi_1}\df\sup_{j\ge0}\left(\sum_{k\ge0}|x_{jk}|^{2}\right)^{\frac12}<\be.
$$

Clearly, $\Xi_p$ is a Banach space for $p\in[1/2,1]$, and $\Xi_p$ is a $(2p)$-Banach space for $p\in(0,\frac12)$.

Let $p\in(0,1]$.
Denote by  $\W_p$ the set of matrices $W$ representable in the form $W=XY^{\rm t}$, where
$X,Y\in\Xi_p$. Put
$$
\|W\|_{\W_p}\df\inf\{\|X\|_{\Xi_p}\|Y\|_{\Xi_p}: ~X,\,Y\in\Xi_p, ~W=XY^{\rm t}\}.
$$
It is well known that $\W_1=\fM_1=\fM_\be$ and $\|\cdot\|_{\W_1}=\|\cdot\|_{\fM_1}=\|\cdot\|_{\fM_\be}$, see, e.g., \cite{Pe}.

Let $q\in(0,\be)$. Clearly, that $\|x+y\|_{\ell^q}^q=\|x\|_{\ell^q}^q+\|y\|_{\ell^q}^q$
for every $x,y\in\ell^{q}$ such that $x\star y=0$ (here we consider vectors $x$ and $y$ as matrices).
This remark implies the following fact.

\begin{lem} 
\label{rqq}
Let $X,Y\in\ell^r(\ell^q)$, where $0<q\le r<\infty$.  Suppose that $X\star Y=0$.
Then $\|X+Y\|_{\ell^r(\ell^q)}^q\le\|X\|_{\ell^r(\ell^q)}^q+\|Y\|_{\ell^r(\ell^q)}^q$.
\end{lem}
\Pf Let $X=\{x_{jk}\}_{j,k\ge0}$ and $Y=\{y_{jk}\}_{j,k\ge0}$. Put $X^{(j)}\df\{x_{jk}\}_{k\ge0}$
and $Y^{(j)}\df\{y_{jk}\}_{k\ge0}$. Then
\bey
\|X+Y\|_{\ell^r(\ell^q)}^q=\left(\sum_{j\ge0}\|X^{(j)}+Y^{(j)}\|_{\ell^q}^r\right)^{q/r}
=\left(\sum_{j\ge0}\Big(\|X^{(j)}\|_{\ell^q}^q+\|Y^{(j)}\|_{\ell^q}^q\Big)^{r/q}\right)^{q/r}\\
\le\left(\sum_{j\ge0}\|X^{(j)}\|_{\ell^q}^r\right)^{q/r}+\left(\sum_{j\ge0}\|Y^{(j)}\|_{\ell^q}^r\right)^{q/r}
=\|X\|_{\ell^r(\ell^q)}^q+\|Y\|_{\ell^r(\ell^q)}^q.\quad\bl
\eey

\begin{thm}
\label{pnorm}
Let $p\in(0,1)$. Then $\W_p$ is a $p$-Banach space. 
\end{thm}

\Pf Let $W=\{w_{jk}\}_{j,k\ge0}$ and $R=\{r_{jk}\}_{j,k\ge0}$ be matrices in $\W_p$. Let us prove that $W+R\in\W_p$ and
$\|W+R\|_{\W_p}^p\le\|W\|_{\W_p}^p+\|R\|_{\W_p}^p$.
Let fix $\e>0$. There exist $X,Y\in\Xi_p$ such that $W=XY^{\rm t}$  and $\|W\|_{\W_p}>\|X\|_{\Xi_p}\|Y\|_{\Xi_p}-\e$. 
Clearly, $w_{jk}=\sum\limits_{l\ge0} x_{jl}y_{kl}$ for all $j,k\ge0$. We can assume that $\|X\|_{\Xi_p}=\|Y\|_{\Xi_p}$.
Moreover, we can assume in addition that $x_{js}=y_{ks}=0$ for all $j,k\ge0$ and all odd $s\ge0$.
In a similar way we can represent $R$ in the form $R=UV^{\rm t}$ in such a way that 
$\|R\|_{\W_p}>\|U\|_{\Xi_p}\|V\|_{\Xi_p}-\e$, $\|U\|_{\Xi_p}=\|V\|_{\Xi_p}$, 
$r_{jk}=\sum\limits_{l\ge0} u_{jl}v_{kl}$ for all $j,k\ge0$ and $u_{jl}=v_{kl}=0$ for all $j,k\ge0$ and all even $l\ge0$.
Clearly, $U\star Y=0$ and $X\star V=0$. Hence, $UY^{\rm t}=0$ and $XV^{\rm t}=0$, whence
$W+R=(X+U)(Y+V)^{\rm t}$. Applying Lemma \ref{rqq} for $q=2p$ and $r=\ps$ we obtain
\bey
\|W+R\|_{\W_p}^{2p}\le\|X+U\|_{\Xi_p}^{2p}\|Y+V\|_{\Xi_p}^{2p}\le(\|X\|_{\Xi_p}^{2p}+\|U\|_{\Xi_p}^{2p})
(\|Y\|_{\Xi_p}^{2p}+\|V\|_{\Xi_p}^{2p})\\
=(\|X\|_{\Xi_p}^{2p}+\|U\|_{\Xi_p}^{2p})^2\le((\|W\|_{\W_p}+\e)^p+(\|R\|_{\W_p}+\e)^p)^2
\eey
Passing to the limit as $\e\to0$, we get $\|W+R\|_{\W_p}^p\le\|W\|_{\W_p}^p+\|R\|_{\W_p}^p$.

It remains to prove that the space $\W_p$ is complete.
It suffices to prove that the series $\sum\limits_{n\ge1}W_n$ converges in $\W_p$ if
$\|W_n\|_{\W_p}<4^{-n}$ for all $n\ge1$. We can take two sequences $\{X_n\}_{n\ge0}$ and $\{Y_n\}_{n\ge0}$
in $\Xi_p$ such that $W_n=X_nY_n^{\rm t}$ and $\|X_n\|_{\Xi_p}=\|Y_n\|_{\Xi_p}<2^{-n}$.
Moreover, it is easy to see that sequences $\{X_n\}_{n\ge0}$ and $\{Y_n\}_{n\ge0}$ can be chosen in
such a way that in addition $X_m\star Y_n=0$ for all $m,n\ge0$ such that $m\ne n$.
Then 
$\sum\limits_{n\ge0}W_n=
\Big(\sum\limits_{n\ge0}X_n\Big)\Big(\sum\limits_{n\ge0}Y_n\Big)^{\rm t}$.
Since the space $\Xi_p$ is complete, it follows that the series 
$\sum\limits_{n\ge0}X_n$ and $\sum\limits_{n\ge0}Y_n$ converge in $\Xi_p$.
Hence, the series $\sum\limits_{n\ge0}W_n$ converges in $\W_p$. $\bl$

\begin{thm}
\label{terrible}
Let $p\in(0,1]$. Then $\W_p\subset\fM_p$ and $\|\cdot\|_{\fM_p}\le\|\cdot\|_{\W_p}$.
\end{thm}

Theorem \ref{terrible} can be reformulated in the following way:

\medskip

{\it Let $p\in(0,1)$ and let $\{a_{jn}\}_{j,n\ge0}$ and $\{b_{nk}\}_{n,k\ge0}$ be matrices
such that the right-hand side of {\em\rf{mune}} is finite. If
$$
\f_{jk}=\sum_{n\ge0}a_{jn}b_{kn},\quad j,\,k\ge0,
$$
and $\Phi=\{\f_{jk}\}_{j,k\ge0}$,
then $\Phi\in\fM_p$ and
\bay
\label{mune}
\|\Phi\|_{\fM_p}
\le\left(\sum_{j\ge0}\left(\sum_{n\ge0}
|a_{jn}|^{2p}\right)^{p_\sharp/p}\right)^{\frac1{2p_\sharp}}
\left(\sum_{k\ge0}\left(\sum_{n\ge0}
|b_{kn}|^{2p}\right)^{p_\sharp/p}\right)^{\frac1{2p_\sharp}}.
\ey
}

\medskip

\Pf It suffices to verify that $\|XY^{\rm t}\|_{\fM_p}\le\|X\|_{\Xi_p}\|Y\|_{\Xi_p}$ for every $X,Y\in\Xi_p$.
In other words we have to prove that 
$$
\|AXY^{\rm t}B\|_{\bS_p}\le\|X\|_{\Xi_p}\|Y\|_{\Xi_p}\|A\|_{\bS_2}\|B\|_{\bS_2}
$$
for every diagonal matrices $A, B\in\bS_2$. Let $\{a_j\}_{j=0}^\be$ and $\{b_k\}_{k=0}^\be$ be the
sequences of diagonal entries of $A$ and $B$. 
Let $X=\{x_{jk}\}_{j,k\ge0}$ and $Y=\{y_{jk}\}_{j,k\ge0}$.
We have

\begin{align*}
\|AXY^{\rm t}B\|^p_{\bS_p}&=
\left\|\left\{\sum_{l\ge0} a_jx_{jl}y_{kl}b_k\right\}_{j,k\ge0}\right\|^p_{\bS_p}
\le\sum_{l\ge0}\Big\|\Big\{a_jx_{jl}y_{kl}b_k\Big\}_{j,k\ge0}\Big\|^p_{\bS_p}\\[.2cm]
&=\sum_{l\ge0}\left(\sum_{j\ge0}|a_j|^2|x_{jl}|^2\right)^{p/2}
\left(\sum_{k\ge0}|b_k|^2|y_{kl}|^2\right)^{p/2}\\[.2cm]
&\le\left(\sum_{l\ge0}\left(\sum_{j\ge0}|a_j|^2|x_{jl}|^2\right)^p\right)^{1/2}
\left(\sum_{l\ge0}\left(\sum_{k\ge0}|b_k|^2|y_{kl}|^2\right)^p\right)^{1/2}\\[.2cm]
&\le\left(\sum_{l\ge0}\sum_{j\ge0}|a_j|^{2p}|x_{jl}|^{2p}\right)^{1/2}
\left(\sum_{l\ge0}\sum_{k\ge0}|b_k|^{2p}|y_{kl}|^{2p}\right)^{1/2}\\[.2cm]
&=\left(\sum_{j\ge0}|a_j|^{2p}\sum_{l\ge0}|x_{jl}|^{2p}\right)^{1/2}
\Big(\sum_{k\ge0}|b_k|^{2p}\sum_{l\ge0}|y_{kl}|^{2p}\Big)^{1/2}\\[.2cm]
&\le\left(\sum_{j\ge0}|a_j|^{2}\right)^{\frac p2}\left(\sum_{j\ge0}
\left(\sum_{l\ge0}|x_{il}|^{2p}\right)^\frac1{1-p}\right)^{\frac{1-p}2}\\[.2cm]
&\times\left(\sum_{k\ge0}|b_k|^{2}\right)^{\frac p2}\left(\sum_{k\ge0}
\left(\sum_{l\ge0}|y_{kl}|^{2p}\right)^\frac1{1-p}\right)^{\frac{1-p}2}\\[.2cm]
&=\|X\|^p_{\Xi_p}\|Y\|^p_{\Xi_p}\|A\|^p_{\bS_2}\|B\|^p_{\bS_2}.
\quad\bl
\end{align*}


\begin{thm} 
\label{terdiag}
Let $p\in(0,1]$ and let $W$ be a diagonal matrix. Then
$\|W\|_{\W_p}=\|W\|_{\fM_p}=\|W\|_{\bS_{\ps}}$.
\end{thm}

\Pf The inequality $\|W\|_{\W_p}\ge\|W\|_{\fM_p}$ follows from Theorem
\ref{terrible}. To prove the opposite inequality we may take the diagonal matrix $X$
such that $X^2=W$, i.e., $W=XX^{\rm t}$. Then 
$\|W\|_{\W_p}\le\|X\|^2_{\Xi_p}=\|X\|^2_{\bS_{2\ps}}=\|W\|_{\bS_{\ps}}=\|W\|_{\fM_p}$ by Theorem \ref{diag}. $\bl$

\begin{thm} 
\label{ter4n}
Let  $W=\{w_{jk}\}_{j,k\ge0}\in\W_p$ with $p\in(0,1)$. Then there exist sequences $\a=\{\a_j\}_{j\ge0}$ and $\b=\{\b_k\}_{k\ge0}$ in $\ell^{2\ps}$ with nonnegative terms
such that $|w_{jk}|\le\a_j\b_k$ for all $j,k\ge0$ and $\|\a\|_{\ell^{2\ps}}\|\b\|_{\ell^{2\ps}}=\|W\|_{\W_p}$.
\end{thm}

\Pf
For each $\e>0$, there exist matrices $X=\{x_{jk}\}_{j,k\ge0}$ and $Y=\{y_{jk}\}_{j,k\ge0}$ 
such that $X,Y\in\Xi_p$, $\|X\|_{\Xi_p}^2=\|Y\|_{\Xi_p}^2<\|W\|_{\W_p}+\e$
and $W=XY^{\rm t}$.
Put 
$$
\a_j(\e)=\left(\sum_{l\ge0}|x_{jl}|^{2p}\right)^{1/2p}\quad
\mbox{and}\quad \b_k(\e)=\left(\sum_{l\ge0}|y_{kl}|^{2p}\right)^{1/2p}.
$$
Then 
$$
\left(\sum_{j\ge 0}(\a_j(\e))^{2\ps}\right)^{1/2\ps}=\left(\sum_{k\ge 0}(\b_k(\e))^{2\ps}\right)^{1/2\ps}<\|W\|_{\W_p}+\e.
$$
We have
$$
|w_{jk}|=\left|\sum_{l\ge0}x_{jl}y_{kl}\right|\le\left(\sum_{l\ge0}|x_{jl}|^2\right)^{1/2}\left(\sum_{l\ge0}|y_{kl}|^2\right)^{1/2}
\le\a_j(\e)\b_k(\e).
$$
Since the sequences $\{\a_j(\e)\}_{j\ge0}$ and $\{\b_k(\e)\}_{k\ge0}$
are uniformly bounded in $\ell^\be$, one can select a sequence $\{\e_\iota\}_{\iota\ge0}$
such that the sequences  $\{\a_j(\e_\iota)\}_{j\ge0}$ and $\{\b_k(\e_\iota)\}_{k\ge0}$
converge in the weak-star topology of $\ell^\be$ to sequences 
$\a=\{\a_j\}_{j\ge0}$ and $\b=\{\b_k\}_{k\ge0}$. It is easy to verify that $\a$ and $\b$ satisfy the requirements. $\bl$

\begin{cor} 
Let $p\in(0,1)$. Then $\W_p\subset\ell^{2\ps}(\Z_+^2)$ and $\|W\|_{\ell^{2\ps}}\le\|W\|_{\W_p}$.
\end{cor}

\begin{cor} 
\label{WLambda}
Let $W=\{w_{jk}\}_{j,k\ge0}$ be a matrix in $\W_p$, where $p\in(0,1]$. Assume that $|w_{jk}|\ge1$
for all $(j,k)\in\L$, where $\L$ is a finite subset $\Z_+^2$. Then
$\|W\|_{\W_p}\ge|\L|^{\frac{1}{2\ps}}$, where $|\L|$ denotes the number of
elements of the set $\L$.
\end{cor}

The following result shows that Theorem \ref{ter4n} is sharp.

\begin{thm} 
\label{ter2n}
Let  $p\in(0,1)$ and let $W=\{w_j\}_{j\ge0}$ be a matrix of rank one. Then
$W\in\W_p$ if and only if $W\in\ell^{2\ps}$. Moreover,
$\|W\|_{\W_p}=\|W\|_{\ell^{2\ps}}$.
\end{thm}

\Pf  Let $w_{jk}=u_jv_k$, $j,k\ge0$, for sequences $\{u_j\}_{j\ge0}$ and $\{v_k\}_{k\ge0}$, and suppose that
$\|W\|_{\ell^{2\ps}}=\|u\|_{\ell^{2\ps}}\|u\|_{\ell^{2\ps}}$.
Let $X=\{u_ja_k\}_{j,k\ge0}$ and $Y=\{v_ja_k\}_{j,k\ge0}$, where $a_0=1$ and $a_j=0$ for all $j\ge1$.
Clearly, $X, Y\in\Xi_p$, $\|X\|_{\Xi_p}=\|u\|_{\ell^{2\ps}}$ and $\|Y\|_{\Xi_p}=\|v\|_{\ell^{2\ps}}$.
We have $XY^{\rm t}=\{u_jv_k\}_{j,k\ge0}$. Hence, $\|W\|_{\W_p}\le\|u\|_{\ell^{2\ps}}\|v\|_{\ell^{2\ps}}=\|W\|_{\ell^{2\ps}}$.
The opposite inequality follows from Theorem \ref{ter4n}. $\bl$

\begin{thm} 
Let $p\in(0,1)$. Then $\ell^\be\otimes_p\ell^\be\not\subset\W_p$ and
$\W_p\not\subset\ell^\be\otimes_p^{\rm c}\ell^\be$.
\end{thm}

\Pf To prove the first statement, it suffices to apply Theorem \ref{ter2n} and observe that under the assumption $\rank W=1$, the matrix
$W$ belongs to $\ell^\be\otimes_p\ell^\be$ if and only if $W\in\ell^\be(\Z^2_+)$. 
The second statement follows from 
Theorems \ref{Dpll} and \ref{terdiag}. $\bl$
\medskip

The next example shows that  in a sense Corollary \ref{WLambda} is sharp.

\medskip

{\bf Example.} Let $W=\{w_{jk}\}_{j,k\ge0}$, where
$w_{jk}=\Bbbone_{S\times T}(j,k)$ for finite subsets $S$ and $T$ of $\Z_+$.
 Then $\|W\|_{\W_p}=|S\times T|^{\frac{1}{2\ps}}$.
By Corollary \ref{WLambda}, we have $\|W\|_{\W_p}\ge|S\times T|^{\frac{1}{2\ps}}$.
To prove the opposite inequality, we can consider the matrices $X=\{x_{jk}\}_{j,k\ge1}$ and
$Y=\{y_{jk}\}_{j,k\ge1}$, where $x_{jk}=\Bbbone_{S\times\{1\}}(j,k)$ and $y_{jk}=\Bbbone_{T\times\{1\}}(j,k)$.
It remains to observe that $W=XY^{\rm t}$ and $\|X\|_{\Xi_p}\|Y\|_{\Xi_p}=|S\times T|^{\frac{1}{2\ps}}$.
Note that $\|W\|_{\fM_p}=\|W\|_{\ell^\be\otimes_p\ell^\be}=1$ if $S\times T\ne\varnothing$.

\

\section{\bf Pisier's sufficient condition}
\setcounter{equation}{0}
\label{Pisier}

\

In this section we analyze another sufficient condition for a matrix to be a Schur multiplier of $\bS_p$ for $p<1$. It is a consequence of a result of Pisier \cite{Pis2}. We compare this sufficient condition with the sufficient condition given in \S\:\ref{uzhas}.

For $q\in(0,\be]$, we define the space $\mY_q$ of scalar matrices by
$$
\mY_q\df\ell^q(\ell^\be)+\big(\ell^q(\ell^\be)\big)^{\rm t}
$$
and we put
$$
\|Z\|_{\mY_q}\df\inf\{\|X\|_{\ell^q(\ell^\be)}+\|Y\|_{\ell^q(\ell^\be)}: Z=X+Y^{\rm t}\},\quad q\in[1,+\be],
$$
$$
\|Z\|_{\mY_q}^q\df\inf\{\|X\|_{\ell^q(\ell^\be)}^q+\|Y\|_{\ell^q(\ell^\be)}^q: Z=X+Y^{\rm t}\},\quad q\in(0,1).
$$

Clearly, $\mY_q$ coincides with the set of matrices $Z=\{z_{jk}\}_{j,k\ge0}$
such that  $|z_{jk}|\le\a_j+\b_k$, $j,k\ge0$, for some nonnegative sequences $\a=\{\a_j\}_{j\ge0}$ and $\b=\{\b_k\}_{k\ge0}$ in $\ell^q$.
We get the same space $\mY_q$ if we require in addition that $\a=\b$.

It is easy to see that $\mY_\be=\ell^\be(\Z^2)$ and $\|\cdot\|_{\mY_\be}=\|\cdot\|_{\ell^\be(\Z^2)}$.

Let us also observe that if
$W$ is a self-adjoint matrix of rank one, 
then $W\in\mY_q$ if and only if $W\in\ell^{2q}(\Z^2)$.
Indeed, suppose that $W=\{u_j\ov u_k\}_{j,k\ge0}$ for a sequence  $\{u_j\}_{j\ge0}$.
Let $\{\a_j\}_{j\ge0}$ be a nonnegative sequence. Then $|u_j\ov u_k|\le\a_j+\a_k$ for all $j,k\ge0$ if and only if $\a_j\ge\frac12|u_j|^2$ for all $j\ge0$.

Recall that in \cite{AP}, with a given number $p$ in $(0,2]$, we associated the number
$$
p_\flat\df\frac{2p}{2-p}.
$$

Let $\fM(\bS_2,\bS_p)$ be the space of matrices $A=\{a_{jk}\}_{j,k\ge0}$ such that
$$
A\star B\in\bS_p \quad\mbox{whenever}\quad B\quad\mbox{is a scalar matrix in}\quad\bS_2.
$$

The following result was obtained by Pisier, see \cite{Pis2}, Theorem 5.1.
$$
\fM(\bS_2,\bS_p)=\mY_{p_\flat}
$$
for all $p\in(0,2]$. 

We need only the easy part of Pisier's result, the inclusion 
\bay
\label{Pisier}
\mY_{p_\flat}\subset\fM(\bS_2,\bS_p),\quad0<p\le2,
\ey
We give here a proof of this inclusion for the reader's convenience.

\medskip

{\bf Proof of \rf{Pisier}.} The result is trivial if $p=2$. Suppose that $p<2$.
Clearly, it suffices to prove that $\|W\star Z\|_{\bS_p}\le\|W\|_{\ell^{p_\flat}(\ell^\be)}\|Z\|_{\bS_2}$
for all $W\in\ell^{p_\flat}(\ell^\be)$ and all $Z\in\bS_2$.
Let $W=\{w_{jk}\}_{j,k\ge0}\in\ell^{p_\flat}(\ell^\be)$. Put $\a_j=\sup_{k\ge0}|w_{jk}|$ and $\a=\{\a_j\}_{j\ge0}$.
Then $\|W\|_{\ell^{p_\flat}(\ell^\be)}=\|\a\|_{\ell^{p_\flat}}$. Let $Z=\{z_{jk}\}_{j,k\ge0}\in\bS_2$.
Using the inequality $\|A\|_{\bS_p}\le\|A\|_{\ell^p(\ell^2)}$ for $p\le2$, we obtain
\begin{align*}
\|W\star Z\|_{\bS_p}^p&\le\sum_{j\ge0}\left(\sum_{k\ge0}|w_{jk}z_{jk}|^2\right)^{\frac p2}
\le\sum_{j\ge0}\a_j^p\left(\sum_{k\ge0}|z_{jk}|^2\right)^{\frac p2}\\[.2cm]
&\le\left(\sum_{j\ge0}\a_j^{\frac{2p}{2-p}}\right)^{\frac{2-p}2}\left(\sum_{j\ge0}\sum_{k\ge0}|z_{jk}|^2\right)^{\frac p2}
=\|W\|_{\ell^{p_\flat}(\ell^\be)}^p\|Z\|_{\bS_2}^p.\quad\bl
\end{align*}

\medskip

It is easy to see that \rf{Pisier} gives us a sufficient condition for a matrix to belong to $\fM_p$.

\begin{cor}
\label{pibemoly}
Let $0<p\le2$. Then $\mY_{p_\flat}\subset\fM_p$.
\end{cor}

Let us prove that this sufficient condition does not cover the sufficient condition 
$\W_p\subset\fM_p$, nor it is covered by the condition $\W_p\subset\fM_p$.

\begin{thm}
\label{nitudainisyuda}
Let $p\in(0,1]$. Then $\W_p\not\subset\mY_{p_\flat}$ and
$\mY_{p_\flat}\not\subset\W_p$.
\end{thm}

\Pf 
Clearly, $\|W\|_{\mY_q}=\|W\|_{\bS_q}$ for any diagonal matrix $W$.
This remark and Theorem \ref{terdiag} implies that $\W_p\not\subset\mY_{p_\flat}$
for $p\in(0,1]$. 

Let us show that $\mY_{p_\flat}\not\subset\W_p$.
Let $W=\{w_{jk}\}_{j.k\ge0}$ with $w_{jk}=\a_j\ne0$ for all $j,k\ge0$, where $\a\in\ell^{p_\flat}$.
Then $W\in\mY_{p_\flat}$. However, $W\not\in\W_p$ by Theorem \ref{ter2n}. $\bl$

\begin{thm} 
\label{thaby}
Let $p\in(0,2]$. Then 
$\mY_q\subset\fM_p$ if and only if $0<q\le p_\flat$.
\end{thm}

\begin{lem} 
\label{pm2}
There exists an $N\times N$ matrix $Z=\{z_{jk}\}_{0\le j,k\le N-1}$
such that $|z_{jk}|=1$ and $\|Z\|_{\fM_p}=N^{\frac1p-\frac12}$ for all $p\in(0,2]$.
\end{lem}

\Pf It is well known that there exists a unitary $N\times N$ matrix $U=\{u_{jk}\}_{0\le j,k\le N-1}$ such
that $|u_{jk}|=N^{-1/2}$. Put $Z=N^{1/2}U$. Denote by $\mathbf{I}_N$ be $N\times N$ matrix with all entries
equal to 1. Then $Z=Z\star\mathbf{I}_N$, whence $\|Z\|_{\bS_p}\le\|Z\|_{\fM_p}\|\mathbf{I}_N\|_{\bS_p}$.
Clearly, $\|Z\|_{\bS_p}=N^{\frac12}\|U\|_{\bS_p}=N^{\frac12+\frac1p}$ and $\|\mathbf{I}_N\|_{\bS_p}=N$. Hence,
$\|Z\|_{\fM_p}\ge N^{\frac1p-\frac12}$. To prove the opposite inequality it suffices to observe that
for every $N\times N$ matrix $A$, we have
$$
\|Z\star A\|_{\bS_p}\le N^{\frac1p-\frac12}\|Z\star A\|_{\bS_2}\le N^{\frac1p-\frac12}\|A\|_{\bS_2}
\le N^{\frac1p-\frac12}\|A\|_{\bS_p}.\quad\bl
$$

\medskip

{\bf Proof of Theorem \ref{thaby}.}  If $q\le p_\flat$, the $\mY_q\subset\mY_{p_\flat}\subset\fM_p$.
It remains to prove that if $\mY_q\subset\fM_p$, then $q\le p_\flat$.

Suppose that $\mY_q\subset\fM_p$ for some $q>p_\flat$. Then there exists a constant $c>0$ such that
$\|W\|_{\fM_p}\le c\|W\|_{\mY_q}$. Let us apply this inequality to the matrix $W=\{w_{jk}\}_{j,k\ge0}$
such that 
$$
w_{jk}\df\left\{\begin{array}{ll}z_{jk},&\text {if}\,\,\,\max(j,k)\le N-1,\\[.2cm]
0,&\text{if}\,\,\,\,\,\,\max(j,k)\ge N,
\end{array}\right.
$$ 
where $Z=\{z_{jk}\}_{0\le j,k\le N-1}$ is the matrix constructed in the proof of Lemma \ref{pm2}.
Then $\|W\|_{\fM_p}=N^{\frac1p-\frac12}$ and $\|W\|_{\mY_q}\le N^{1/q}$. Hence,
$N^{\frac1p-\frac12}\le c N^{1/q}$ for all $N\ge1$, and so $p_\flat^{-1}\le q^{-1}$, i.e. $q\le p_\flat$. $\bl$

\medskip

{\bf Remark 1.} It is easy to see that $\ell^\be\otimes_p\ell^\be\not\subset\mY_{p_\flat}$
for $p\in(0,1]$. Indeed, the infinite matrix $\mathbf{I}_\be$ with entries identically equal to 1
belongs to $\ell^\be\otimes_p\ell^\be$ and does not  belong to $\mY_{p_\flat}$.
In the case $p\in(0,1)$, the authors do not know whether $\mY_{p_\flat}\subset\ell^\be\otimes_p\ell^\be$ or not.

\medskip

We conclude this section with a few observations.

Let $A=\{a_{jk}\}_{j,k\ge0}$ be a Schur multiplier of $\bS_p$. It is said to be
an {\it absolute Schur multiplier of} $\bS_p$ if every matrix $B=\{b_{jk}\}_{j,k\ge0}$ satisfying $|b_{jk}|\le|a_{jk}|$, $j,k\ge0$, is also a Schur multiplier of $\bS_p$.
Denote by $(\fM_p)_{\rm C}$ the {\it space of absolute Schur multipliers of} $\bS_p$.

We can use a similar same notation in a more general situation.
Let $\mX$ be linear subset of $\ell^\be(\Z_+^2)$. Put 
$$
\mX_{\rm C}\df\{X\in \mX: X\star\ell^\be(\Z_+^2)\subset \mX\}.
$$

Clearly, $A\in(\fM_p)_{\rm C}$ if and only if $A$ is a Schur multiplier from $\bS_p$ into $(\bS_p)_{\rm C}$.
It is well known that $(\bS_p)_{\rm C}=\ell^p(\ell^2)+\big(\ell^p(\ell^2)\big)^{\rm t}$ for $p\in(0,2]$, see \cite{Pis2} in the case $p<1$, the cases $p\in(1,2)$ and $p=1$ were considered in \cite{LP} and \cite{LPP}. 

It is easy to see that for $p\in(0,2]$, every matrix $A$ in $\fM(\bS_2,\bS_p)$ is an absolute Schur multiplier of $\bS_p$,
i.e. $\fM(\bS_2,\bS_p)\subset(\fM_p)_{\rm C}$ for $p\in(0,2]$.
Hence, in Corollary \ref{pibemoly}
and Theorem \ref{thaby} the space $\fM_p$ can be replaced with the space $(\fM_p)_{\rm C}$.

\medskip

{\bf Remark 2.}  Let $p\in(0,2)$. Then $\fM(\bS_2,\bS_p)\ne(\fM_p)_{\rm C}$.
Indeed, it is easy to see that a diagonal matrix belongs to $\fM_p$ if and only if
it belongs to
$(\fM_p)_{\rm C}$. 
It remains to observe that a diagonal matrix $W$ belongs to $\fM(\bS_2,\bS_p)$ if and only if
$W\in\bS_{p_\flat}$ and apply Theorem \ref{diag}.

\medskip

This remark implies 
$\mY_{p_\flat}\ne(\fM_p)_{\rm C}$ for $p\in(0,2)$ which complements Corollary \ref{pibemoly}.

\medskip

{\bf Remark 3.}  Let $p\in(0,2)$.  Lemma \ref{pm2} allows us to prove that
there exist a matrix $A=\{a_{jk}\}_{j,k\ge0}$ and two sequences
$\a=\{\a_j\}_{j\ge0}$ and $\b=\{\b_k\}_{k\ge0}$ in $\ell^{2\ps}$
with nonnegative terms such that  $|a_{jk}|\le\a_j\b_k$ for all $j,k\ge0$
but $A\not\in\W_p$, see the proof of Theorem \ref{thaby}.

\medskip

This remark implies that $\W_p\not\subset(\fM_p)_{\rm C}$ for $p\in(0,1)$.
This is another way to see that $\W_p\not\subset\mY_{p_\flat}$ for $p\in(0,1)$, see Theorem \ref{nitudainisyuda}.

\

\
 
\noindent
\begin{tabular}{p{7cm}p{15cm}}
A.B. Aleksandrov & V.V. Peller \\
St.Petersburg Branch & Department of Mathematics \\
Steklov Institute of Mathematics  & Michigan State University \\
Fontanka 27, 191023 St.Petersburg & East Lansing, Michigan 48824\\
Russia&USA\\
email: alex@pdmi.ras.ru&and\\
&Peoples' Friendship University\\
& of Russia (RUDN University)\\
&6 Miklukho-Maklaya St., Moscow,\\
& 117198, Russian Federation\\
& email: peller@math.msu.edu
\end{tabular}

\end{document}